\documentclass[conference]{IEEEtran}

\usepackage{graphicx,times,amsmath} 
\usepackage{epsf,epsfig}
\usepackage{graphics}
\usepackage{subfigure}
\usepackage{multicol}
\usepackage{amsmath,amssymb}

\def\R{{\mathbb R}}

\IEEEoverridecommandlockouts	
				
\newcommand{\vecc}[1]{\mathbf{#1}}

\newcommand{\xmean}[1]{\langle\vec{x}\rangle_{\textrm{W}}^{(#1)}}
\newcommand{\mat}[1]{\mathchoice{\mbox{\boldmath$\displaystyle#1$}}
  {\mbox{\boldmath$\textstyle#1$}} {\mbox{\boldmath$\scriptstyle#1$}}
  {\mbox{\boldmath$\scriptscriptstyle#1$}}}

\newcommand{\cc}{\vecc{c}}
\newcommand{\RR}{\mathbb{R}}
\newcommand{\HH}{\vecc{H}}
\newcommand{\FF}{\vecc{F}}
\newcommand{\NN}{\vecc{N}}
\newcommand{\KK}{\vecc{K}}
\newcommand{\EE}{\vecc{E}}
\newcommand{\cost}{\vecc{\mathcal{J}}}

\newcommand{\z}{\hspace{0.11785em}} 

\newcommand{\eq}{Eq.\z\z\ref}

\begin{document}

\title{\ \\ \LARGE\bf Identification of the Isotherm Function in Chromatography \\
 Using CMA-ES
 \thanks{1. TAO Project-Team -- INRIA Futurs, LRI, Orsay, 2. Math\'ematiques, Applications et Physique Math\'ematique d'Orl\'eans, 3. Laboratoire Jacques-Louis Lions -- UPMC, Paris, {\tt  \{jebalia,auger,marc\}@lri.fr}, {\tt Francois.James@math.cnrs.fr}, {\tt postel@ann.jussieu.fr} 
 }
}

\author{M. Jebalia$^{1}$, A. Auger$^{1}$, M. Schoenauer$^{1}$,  F. James$^{ 2}$, M.Postel$^{3}$}

\maketitle

\begin{abstract}
This paper deals with the identification of the flux for a system of
conservation laws in the specific example of analytic chromatography.
The fundamental equations of chromatographic process are highly
non linear. The state-of-the-art Evolution Strategy, CMA-ES (the
Covariance Matrix Adaptation Evolution Strategy), is used to identify
the parameters of the so-called isotherm function. The approach was validated on different
configurations of simulated data using either one, two or three
components mixtures. CMA-ES is then applied to real data cases and its
results are compared to those of a gradient-based strategy.
\end{abstract}


\section{Introduction}
The chromatography process is a powerful tool to separate or analyze
mixtures \cite{G06}. It is widely used in chemical industry (pharmaceutical,
perfume and oil industry, etc) to produce relatively high
quantities of very pure components. This is achieved by taking advantage of the selective absorption of the different components in a solid porous medium. The moving fluid mixture is
percolated through the motionless medium in a column. The various
components of the mixture propagate in the column at different
speeds, because of their different affinities with the solid
medium. The art of chromatography separation requires predicting the
different proportions of every component of the mixture at the end of
the column (called {\it the chromatogram}) during the experiment.  In
the ideal (linear) case, every component has its own fixed propagation
speed, that does not depend on the other components. In this case, if
the column is sufficiently long, pure components come out at the end
of the column at different times: they are perfectly separated. But in
the real world, the speed of a component heavily depends on every
other component in the mixture. Hence, the fundamental Partial
Differential Equations of the chromatographic process, derived from
the mass balance, are highly non linear.  The process is governed by
a nonlinear function of the mixture concentrations, the so-called {\it
Isotherm Function}. This function computes the amount of absorbed
quantity of each component w.r.t. all other components.

Mathematically speaking, thermodynamical properties of the isotherm
ensure that the resulting system of PDEs is hyperbolic, and standard
numerical tools for hyperbolic systems can hence be applied; if
the isotherm is known: The precise knowledge of the isotherm is
crucial, both from the theoretical viewpoint of physico-chemical modeling and 
regarding the more practical preoccupation of accurately controlling 
the experiment to improve separation.  Specific chromatographic
techniques can be used to directly identify the isotherm, but gathering a few
points requires several months of careful experiments. Another
possible approach to isotherm identification consists in solving the
inverse  problem numerically: find the isotherm such that numerical
simulations result in chromatograms that are as close as possible
to the actual experimental outputs.

This paper introduces an evolutionary method to tackle the
identification of the isotherm function from experimental
chromatograms. The goal of the identification is to minimize the
difference between the actual experimental chromatogram and the
chromatogram that results from the numerical simulation of the
chromatographic process. Chemical scientists have introduced several
parametric models for isotherm functions (see \cite{G06} for all
details of the most important models). The resulting optimization
problem hence amounts to parametric optimization, that is addressed
here using the state-of-the-art Evolution Strategy,
CMA-ES. Section~\ref{Chromatography} introduces the direct problem and
Section~\ref{Optimization} the optimization (or inverse) problem. Section
\ref{State-of-the-art} reviews previous approaches to the problem based on gradient
optimization algorithms \cite{cit:james,James:2007}. Section \ref{CMA} details the CMA-ES
method and the implementation used here. Finally,
Section~\ref{sec:results} presents experimental results: first,
simulated data are used to validate the proposed approach; second, real
data are used to compare the evolutionary approach with a
gradient-based method.

\section{Physical problem and model}
\label{Chromatography}

Chromatography aims at separating the components of a mixture based on the selective absorption of chemical species by a solid porous 
medium. The  fluid mixture moves down  through a column of length $L$, considered here to be one-dimensional. The various components of the mixture propagate in the column at different speeds, because of their different behavior
when interacting   with the porous medium. 
At a given time $t \in \R^+$, for a given $z \in [0,L]$ the concentration of $m$ species is a real vector of $\R^m$ 
denoted $\cc(t,z)$. The evolution of $\cc$ is governed by the following partial differential equation:
 \begin{equation}
 \left\{
 \begin{aligned}
 \partial_z \cc + \partial_t \FF(\cc) & =  0, \\
     \cc(0,z) & =  \cc_{0}(z),     \\
 \cc(t,0) & =  \cc_{inj}(t).
 \label{eq:PDEc}
 \end{aligned}
 \right.
 \end{equation}
where $\cc_0: \RR \rightarrow \RR^m$ is the initial concentration, $\cc_{inj}: \RR \rightarrow \RR^m$ the injected concentration at the entrance of the column and $\FF: \RR^m \rightarrow \RR^m$ is the flux function that can be expressed in the following way
$$
\FF(\cc) = \frac1u \left(\cc + \frac{1-\epsilon}{\epsilon} \HH(\cc) \right)
$$
where $\HH: \RR^m \rightarrow \RR^m$ is the so-called isotherm
function, $\epsilon \in (0,1)$ and $u \in \RR^+$
\cite{James:2007}. The Jacobian matrix of $\FF$ being diagonalizable
with strictly positive eigenvalues, the system (\ref{eq:PDEc}) is
strictly hyperbolic and thus admits an unique solution as soon as
$\FF$ is continuously differentiable, and the initial and
injection conditions are piecewise continuous. The solution of
\eq{eq:PDEc} can be approximated using any finite difference method
that is suitable for hyperbolic systems \cite{Raviart91}. A uniform
grid in space and time of size $(K+1) \times (N+1)$ is defined: Let
$\Delta z$ (resp. $\Delta t$) such that $K \Delta z = L$ (resp. $N
\Delta t = T$).  Then an approximation of the solution of \eq{eq:PDEc}
can be computed with the Godunov scheme:
\begin{equation}
\cc_{k+1}^{n}=\cc_{k}^{n} -  \frac{\Delta z}{\Delta t}(\FF(\cc_{k}^{n})-\FF(\cc_{k}^{n-1})) 
\label{eq:godunov-scheme}
\end{equation}
where $\cc_k^n $ is an approximation of the mean value of the solution
$\cc$ at point $(k \Delta z, n \Delta t)$\footnote{Mean value over the
volume defined by the corresponding cell of the grid.}. For a fixed
value of $\frac{\Delta z}{\Delta t}$, the solution of
\eq{eq:godunov-scheme} converges to the solution of \eq{eq:PDEc} as
$\Delta t$ and $\Delta z$ converge to zero. The numerical scheme given
in \eq{eq:godunov-scheme} is numerically stable under the so-called
CFL condition stating that the largest absolute value of the
eigenvalues of the Jacobian matrix of $\FF$ is upper-bounded by a
constant
\begin{equation}
\frac{\Delta z}{\Delta t} \max_c \textrm{Sp} (|\FF'(c)|) \leq \textrm{CFL} < 1.
\label{eq:cfl}
\end{equation}

\section{The Optimization Problem}
\label{Optimization}
\subsection{Goal}
\label{sec:goal}
The goal is to identify the isotherm function from experimental
chromatograms: given initial data $c_0$, injection data $c_{inj}$, and
the corresponding experimental chromatogram $c_{exp}$ (that can be
either the result of a simulation using a known isotherm function, or
the result of actual experiments by chemical scientists), find the
isotherm function $\HH$ such that the numerical solution of  \eq{eq:PDEc}
using the same initial and injection
conditions results in a chromatogram as close as possible to the
experimental one $c_{exp}$.

Ideally, the goal is to find $\HH$ such that the following system of
PDEs has a unique solution $\cc(t,z)$:
 \begin{equation}
 \left\{
 \begin{aligned}
 \partial_z \cc + \partial_t \FF(\cc) & =  0, \\
     \cc(0,z) & =  \cc_{0}(z),     \\
 \cc(t,0) & =  \cc_{inj}(t),\\
 \cc(t,L) & =  \cc_{exp}(t).
 \label{eq:PDEinverse}
 \end{aligned}
 \right.
 \end{equation}
However, because in most real-world cases this system will not have an
exact solution, it is turned into a minimization problem. For a given
isotherm function $\HH$, solve system \ref{eq:PDEc} and define the
cost function $\cost$ as the least square difference between the
computed chromatogram $\cc_{\HH}(t,L)$ and the experimental
one $\cc_{exp}(t)$:
\begin{equation}
\label{eq:fitness}
\cost(\HH)=\int_{0}^T \| \cc_{\HH}(t,L) - \cc_{exp}(t) \|^2 dt
\end{equation}
%
If many experimental chromatograms are provided, the cost function is 
the sum of such functions   $\cost$  computed for each experimental chromatogram. 
\subsection{Search Space}
When tackling a function identification problem, the first issue to
address is the parametric vs non-parametric choice \cite{Eurogen2001}:
parametric models for the target function result in parametric
optimization problems that are generally easier to tackle -- but a bad
choice of the model can hinder the optimization. On the other hand,
non-parametric models are a priori less biased, but search algorithms
are also less efficient on large unstructured search space.

Early trials to solve the chromatography inverse problem using a
non-parametric model (recurrent neural-network) have brought a
proof-of-concept to such approach \cite{cit:fadamarc}, but have also
demonstrated its limits: only limited precision could be reached, and
the approach poorly scaled up with the number of components of the
mixture.
 
Fortunately, chemists provide a whole zoology of parametrized models
for the isotherm function $\HH$, and using such models, the
identification problem amounts to parametric optimization. For $i
\in \{1, \ldots, m \}$, denote $\HH_i$ the component $i$ of the
function $\HH$. The main models for the isotherm function that will be
used here are the following:
\begin{itemize}
\item 	The {\bf Langmuir} isotherm \cite{Langmuir1918} assumes that the different components are in competition to occupy each site of the porous medium. This gives, for all $i=1, \ldots, m$ 
\begin{equation}
\label{eq:lang}
\HH_i(c)=\frac{\NN^*}{1+\sum_{l=1}^m  \KK_l \cc_l } \ \KK_i  \cc_i.
\end{equation}
There are $m+1$ positive parameters: the {\it Langmuir coefficients}
$(\KK_i)_{i \in [1, m]}$, homogeneous to the inverse of a
concentration, and the {\it saturation coefficient} $\NN^*$ that
corresponds to some limit concentration.

\item   The {\bf Bi-Langmuir} isotherm generalizes the Langmuir isotherm by assuming two different kinds of sites on the absorbing medium. The resulting equations are, for all $i=1, \ldots, m$ 
\begin{equation}
\label{eq:bi-lang}
 \HH_i(c)=\sum_{s \in \{1,2\}}  \frac {\NN^*_s}{1+\sum_{l=1}^m \ \KK_{l,s} \cc_l}   \ \KK_{i,s} \cc_i.
\end{equation}
This isotherm function here depends on $2(m+1)$ parameters: the
generalized Langmuir coefficients ${(\KK_{i,s})}_{i \in [1, m],
s=1,2}$ and the generalized saturation coefficients
$(\NN^*_s)_{s=1,2}$.

\item  The  {\bf Lattice} isotherm \cite{valentinJames97} is a generalization of Langmuir isotherm that also considers interactions among the different sites of the porous medium. Depending on the degree $d$ of interactions (number of interacting sites grouped together), this model depends, additionally to the Langmuir coefficients $(\KK_i)_{i \in [1,m]}$ and the saturation coefficient $\vecc{N}^*$, on interaction energies $(\EE_{ij})_{i,j \in  [0, d], 2 \leq i+j \leq d}$ resulting in  $\prod_{i=1}^m  \frac{d+i}{i} $ parameters. For instance, for  one component ($m=1$) and degree $2$, this gives:
\begin{equation}
\label{eq:lattice}
\HH_1(\cc)=\frac{\vecc{N}^*}{2} \frac{ \KK_1 \  \cc + e^{-\frac{\EE_{11}}{RT}}( \KK_1 \  \cc)^2}{1+ 2  \KK_1 \  \cc + e^{-\frac{\EE_{11}}{RT}}( \KK_1 \  \cc)^2},
\end{equation}
where $T$ is the absolute temperature and $R$ is the universal gas
constant. Note that in all cases, a Lattice isotherm with $0$ energies
simplifies to the Langmuir isotherm with the same Langmuir and
saturation coefficients up to a factor $\frac{1}{2}$.

\end{itemize}


\section{Approach Description}
\subsection{Motivations}
\label{State-of-the-art}

Previous works on parametric optimization of the chromatography inverse problem have used gradient-based approaches
\cite{cit:james,James:2007}. In \cite{cit:james}, the
gradient of $\cost$ is obtained by writing and solving numerically the
adjoint problem, while direct differentiation of the discretized equation 
have also been investigated in \cite{James:2007}. However the fitness function to optimize
is not necessarily convex and no results are provided for
differentiability. Moreover,  experiments performed in 
\cite{James:2007} suggest that the function is multimodal, since the 
gradient algorithm converges to different local optima depending on
the starting point. Evolutionary algorithms (EAs) are stochastic
global optimization algorithms, less prone to get stuck in local
optima than gradient methods, and do not rely on convexity
assumptions. Thus they seem a good choice to tackle this problem.  Among
EAs, Evolution Strategies have been specifically designed for
continuous optimization. The next section introduces the state of the art
EA for continuous optimization, the covariance matrix adaptation ES
(CMA-ES).


\subsection{The CMA Evolution Strategy}
\label{CMA}

CMA-ES is a stochastic optimization algorithm specifically designed
for continuous optimization
\cite{Hansen:2001,Hansen:2003,Hansen:2004b,Auger:2005a}. At each iteration $g$, a population of points of an $n$-dimensional
continuous search space (subset of $\RR^n$), is sampled according to a
multi-variate normal distribution. Evaluation of the fitness of the
different points is then performed, and parameters of the multi-variate
normal distribution are updated. 

More precisely, let $\xmean{g}$ denotes the mean value of the 
(normally) sampling 
distribution at iteration $g$. Its covariance matrix is usually factorized 
in two terms: 
$\sigma^{(g)} \in \mathbb{R}^+$, also called the {\it step-size}, and $\mat{C}^{(g)}$, a definite positive $n \times n$ matrix,  that
is abusively called the covariance matrix. The independent sampling of the $\lambda$ offspring can then be written:
\begin{equation*}
\label{eq:offspring}
  \vec{x}_k^{(g+1)} = \xmean{g} + \mathcal{N}_k \left(0, (\sigma^{(g)})^2
      \mat{C}^{(g)} \right) \mbox{ for } k=1,\ldots,\lambda
\end{equation*}
where $ \mathcal{N}_k \left(0, M \right)$ denote independent realizations of the multi-variate normal distribution of covariance matrix $M$.

The $\mu$ best offspring are recombined into
\begin{equation}
\label{eq:updatemean}
\xmean{g+1} = \sum_{i=1}^\mu
w_i \vec{x}_{i:\lambda}^{(g+1)} \enspace,
\end{equation}
where the positive weights $w_i\in\mathbb{R}$ are set according to individual ranks and sum to one. The index $i\!:\!\lambda$ denotes the $i$-th best offspring. \eq{eq:updatemean} can be rewritten as
\begin{equation}
\label{eq:updatemean2}
\xmean{g+1} = \xmean{g} + \sum_{i=1}^\mu w_i \mathcal{N}_{i:\lambda} \left(0, (\sigma^{(g)})^2
      \mat{C}^{(g)} \right) \enspace,
\end{equation}
The covariance matrix $\mat{C}^{(g)}$ is a positive definite symmetric matrix. Therefore it can be decomposed in
$$
\mat{C}^{(g)} = \mat{B}^{(g)}\mat{D}^{(g)} \mat{D}^{(g)} \left(\mat{B}^{(g)}\right)^T \enspace,
$$
where $\mat{B}^{(g)}$ is an orthogonal matrix, {\it i.e.} $\mat{B}^{(g)} \left(\mat{B}^{(g)}\right)^T = I_d$ and $\mat{D}^{(g)}$ a diagonal matrix whose diagonal contains the square root of the eigenvalues of $\mat{C}^{(g)}$.

The so-called strategy parameters of the algorithm, the covariance matrix
$\mat{C}^{(g)}$ and the step-size $\sigma^{(g)}$, are updated so as to
increase the probability to reproduce good steps. The so-called
rank-one update for $\mat{C}^{(g)}$ \cite{Hansen:2001} takes place as follows. First, an
evolution path is computed:
\begin{equation*}
 \vec{p}_c^{(g+1)} =(1- c_{c}) \vec{p}_c^{(g)} + \frac{\sqrt{c_{c}(2-c_{c})\mu_{\rm eff}}}{\sigma^{(g)}} \left(\xmean{g+1} - \xmean{g} \right)
\end{equation*}
where $c_c \in ]0,1]$ is the cumulation coefficient and $\mu_{\rm
eff}$ is a strictly positive coefficient. This evolution path can be
seen as the descent direction for the algorithm. 

Second the covariance
matrix $\mat{C}^{(g)}$ is ``elongated`` in the direction of the
evolution path, {\it i.e. } the rank-one matrix $\vec{p}_c^{(g+1)}
\left(\vec{p}_c^{(g+1)}\right)^{T}$ is added to $\mat{C}^{(g)}$:
\begin{align*}
\mat{C}^{(g+1)} & =(1-c_{\rm cov}) \mat{C}^{(g)} +  c_{\rm cov} \vec{p}_c^{(g+1)} \left(\vec{p}_c^{(g+1)}\right)^{T} \enspace
\end{align*}
where $c_{\rm cov} \in ]0,1[$. The complete update rule for the covariance matrix is a combination of
the rank-one update previously described and the rank-mu update
presented in \cite{Hansen:2003}. 

The update rule for the step-size
$\sigma^{(g)}$ is called the path length control. First, another evolution path is computed:
\begin{multline}
 \vec{p}_\sigma^{(g+1)} =(1- c_{\sigma}) \vec{p}_\sigma^{(g)} + \frac{\sqrt{c_\sigma(2-c_\sigma)\mu_{\rm eff}}}{\sigma^{(g)}} \times \\ \mat{B}^{(g)}{\mat{D}^{(g)}}^{-1}{\mat{B}^{(g)}}^{T} \left(\xmean{g+1} - \xmean{g} \right)
\end{multline}
where $c_\sigma \in ]0,1]$.  The length of this vector is compared to
the length that this vector would have had  under random selection,
{\it i.e.} in a scenario where no information is gained from the
fitness function and one is willing to keep the step-size
constant. Under random selection the vector $\vec{p}_\sigma^{(g)}$ is
distributed as $\mathcal{N}(0,I_d)$. Therefore, the step-size is
increased if the length of $\vec{p}_\sigma^{(g)}$ is larger than
$\EE(\parallel \mathcal{N}(0,I_d)\parallel)$ and decreased if it is
shorter. Formally, the update rule reads:
\begin{equation}
\sigma^{(g+1)}=\sigma^{(g)}
\exp \left(
\frac{c_\sigma}{d_{\sigma}}
\left( 
\frac{\parallel \vec{p}_{\sigma}^{(g+1)}\parallel }
{\EE(\parallel \mathcal{N}(0,I_d)\parallel)} 
-1\right)\right)
\label{eq:eq34}
\end{equation}
where $d_\sigma > 0 $ is a damping factor.

The default parameters for CMA-ES were carefully derived in \cite{Hansen:2004b}, Eqs.~6-8. The only problem-dependent parameters are $\xmean{0}$ and $\sigma^{(0)}$, and, to some extend, the offspring size $\lambda$: its default value is 
$\lfloor 4 + 3 \log(n) \rfloor$ (the $\mu$ default value is $\lfloor  \frac{\lambda}{2} \rfloor$), but increasing  $\lambda$
increases the probability to converge towards the global optimum when
minimizing multimodal fitness functions \cite{Hansen:2004b}. 

This fact was systematically exploited in \cite{Auger:2005a}, 
where a "CMA-ES restart" algorithm is proposed, in which the
population size is increased after each restart. Different
restart criteria are used:
\begin{enumerate}
\item {\it RestartTolFun}: Stop if the range of the best objective function values of the recent generation is below than 
a TolFun value.
\item {\it RestartTolX}: Stop if the standard deviation of the normal distribution is smaller than a TolX value and $\sigma 
\vec{p}_c$ is smaller than TolX in all components.
\item {\it RestartOnNoEffectAxis}: Stop if adding a $0.1$ standard deviation vector in a principal axis direction of $\mat{C}^
{(g)}$ does not change $\xmean{g}$.
\item {\it RestartCondCov}: Stop if the condition number of the covariance matrix exceeds a fixed value.
\end{enumerate}
The resulting algorithm (the CMA-ES restart, simply denoted CMA-ES in the remainder of this paper) is a quasi parameter free algorithm
that performed best for the CEC 2005 special session on parametric
optimization \cite{CMARestart}.

An important 
property of CMA-ES is its invariance to linear transformations of the
search space. Moreover, because of the rank-based selection, CMA-ES is
invariant to any monotonous transformation of the fitness function: 
optimizing $f$ or $h \circ f $ is equivalent, for any rank-preserving function $h:\R \rightarrow \R$. In particular, convexity has no impact on the actual  behavior of CMA-ES.


\subsection{CMA-ES Implementation}
This section describes the specific implementation of CMA-ES to identify $n$ isotherm coefficients. 
For the sake of clarity we will
use a single index in the definition of the coefficients of the isotherm,
{\it i.e} we will identify $\KK_{a}$, $\NN^{*}_{b}$ and $\EE_{c}$ for 
$a \in [1, A]$, $b \in
[1, B]$ and $c \in [1, C]$ where $A$, $B$ and $C$
are integers summing up to $n$. \\

\paragraph*{\bf Fitness function and CFL condition} 
The goal is to minimize the fitness function defined in Section  \ref{sec:goal}. 
In the case where identification is done using only one experimental chromatogram, 
the fitness function is the function $\cost$ defined in
\eq{eq:fitness} as the least squared difference between an
experimental chromatogram $\cc_{exp}(t)$ obtained using experimental conditions $\cc_0$, $\cc_{inj}$ and a numerical approximation
of the solution of system (\ref{eq:PDEc}) for a candidate isotherm
function $\HH$ using the same experimental conditions. The numerical simulation of a
solution of \eq{eq:PDEc} is computed with a Godunov scheme written
in C++ (see \cite{CHROMALGEMA} for the details of the implementation). 

In order
to validate the CMA-ES approach, first "experimental" chromatograms  were in fact
computed using numerical simulations of
\eq{eq:PDEc} with different experimental conditions. Let $\FF_{sim}$ denotes
the flux function used to simulate the experimental chromatogram. For
the simulation of an approximated solution of \eq{eq:PDEc}, a time
step ${\Delta t}$ and a CFL coefficient strictly smaller than
one (typically 0.8) are fixed beforehand.  The quantity $\textrm{max Sp} (|\FF_{sim}'(c)|)$
is then estimated using a power method, and the space step
${\Delta z}$ can then be set such that \eq{eq:cfl} is satisfied for $\FF_{sim}$. The same ${\Delta t}$
and ${\Delta z}$ are then used during the optimization of $\cost$.

When $\cc_{exp}$ comes from real data, an initial value for the parameters to estimate, i.e. an initial	 guess given by the
expert is used to set the CFL condition (\ref{eq:cfl}).\\



\begin{figure}
\centerline{\subfigure[Simulation using a Langmuir isotherm, identification using a Langmuir model: the chromatogram is perfectly fit.]
{\includegraphics[width=1.5in]{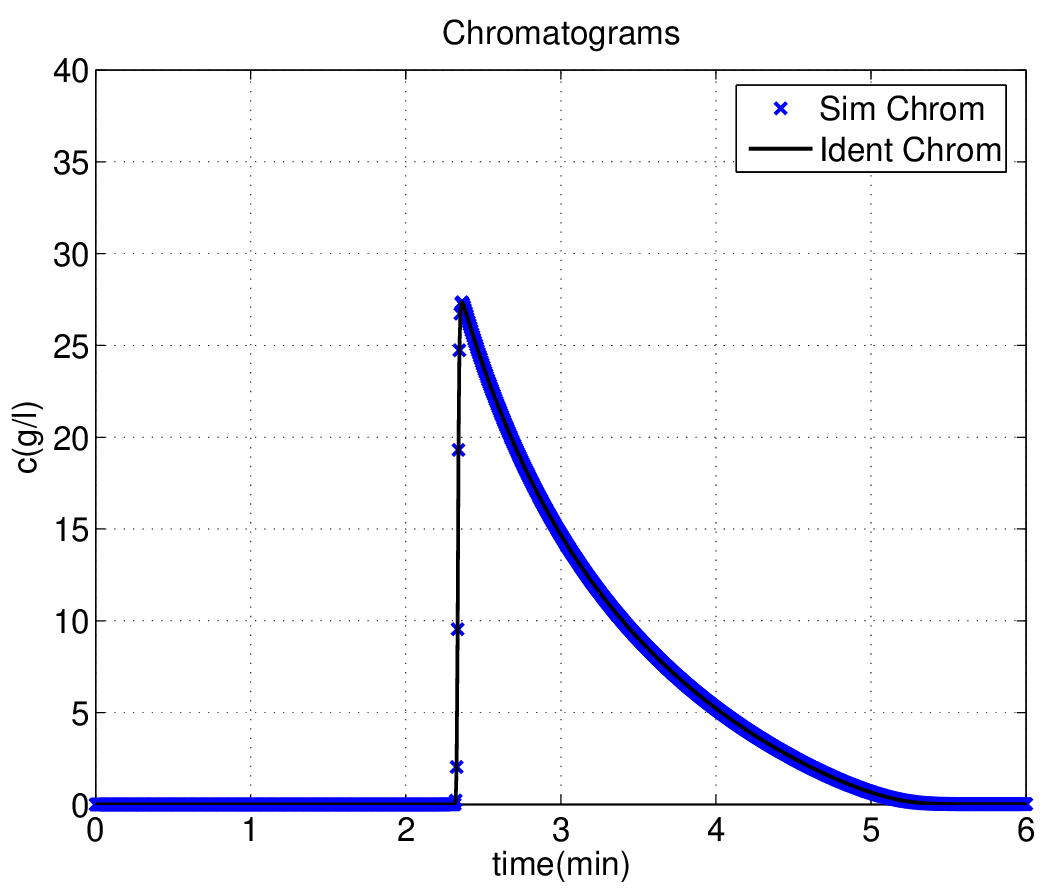}
\label{fig_sim1}}
\hfil
\subfigure[Simulation using a Lattice isotherm, identification using a Langmuir model: poor fit of the chromatogram.]{\includegraphics[width=1.5in]
{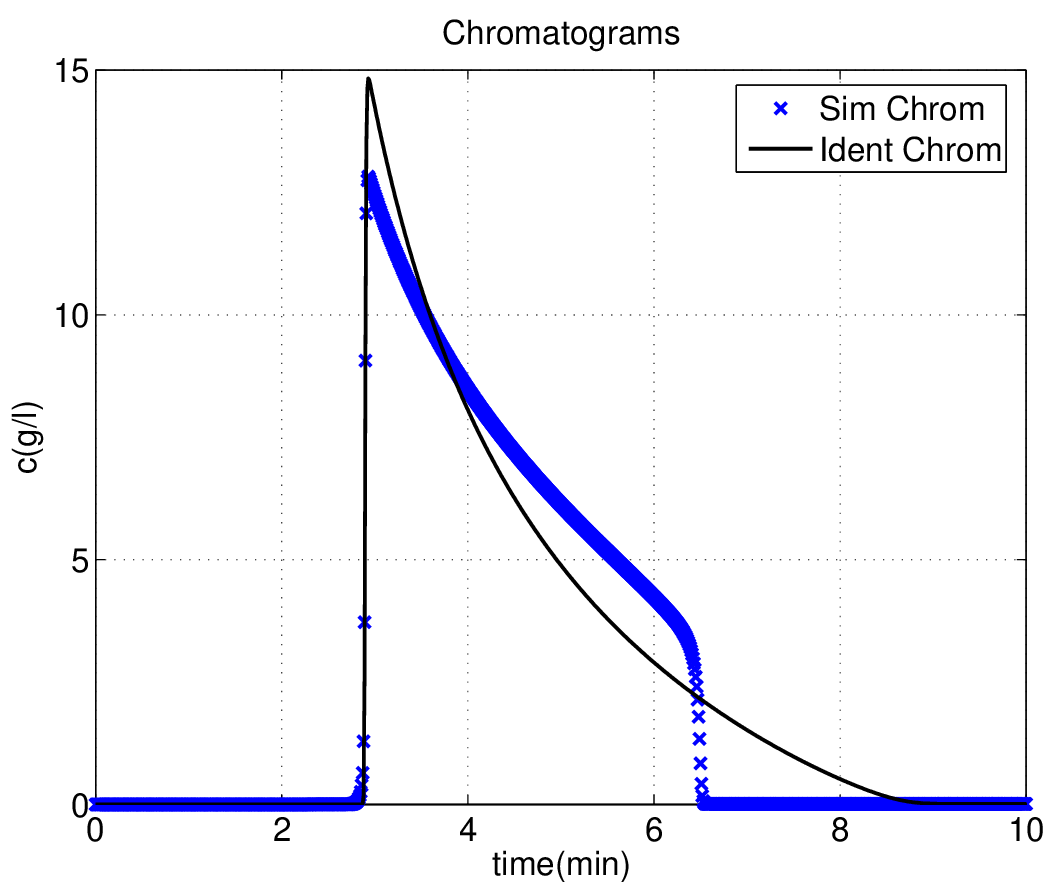}
\label{fig_sim2}}}
\caption{Importance of the choice of model (one component mixture)}
\label{fig_sim}
\end{figure}

\paragraph*{\bf Using expert knowledge} 
The choice of the type of isotherm function to be identified 
will be, in most cases, given by the chemists. Fig~\ref{fig_sim}
illustrates the importance of this choice. In
Fig~\ref{fig_sim}-(a), the target chromatogram $\cc_{exp}$ is
computed using a Langmuir isotherm with one component
($m=1$ and thus $n=2$). In Fig~\ref{fig_sim}-(b), the target chromatogram $\cc_{exp}$
is computed using a Lattice of degree $3$ with one component ($m=1$ and thus $n=4$). In both cases, the
identification is done using a Langmuir model, with $n=2$. 
It is clear from the figure that one is able to correctly identify the isotherm, and hence fit the
"experimental" chromatogram when choosing the correct model (Fig~\ref{fig_sim} (a)) whereas the fit of the chromatogram is very poor when the
model is not correct (Fig~\ref{fig_sim} (b)).


Another important issue when using CMA-ES is the initial choice for the
covariance matrix: without any information, the algorithm starts with
the identity matrix. However, this is a poor choice in case the
different variables have very different possible order of magnitude, and the
algorithm will spend some time adjusting its principal directions to
those ranges.

In most cases of chromatographic identification, however, chemists provide orders of magnitudes, bounds and
initial guesses for the different values of the unknown parameters. 
Let $[(\KK_{a})_{min},(\KK_{a})_{max}]$,
$[(\NN_{b}^{*})_{min},(\NN_{b}^{*})_{max}]$ and
$[(\EE_{c})_{min},(\EE_{c})_{max}]$ the ranges guessed by the chemists
for respectively each $\KK_a$, $\NN^{*}_b$ and $\EE_c$. All parameters 
are linearly scaled into those intervals from $[-1, 1]$, removing the need to 
modify the initial covariance matrix of CMA-ES.\\

\paragraph*{\bf Unfeasible solutions} 
Two different situations can lead to {\it unfeasible} solutions:

First when one parameter at least, among parameters which have to  be positive,  becomes negative (remember that CMA-ES generates offspring using an unbounded normal distribution),  the fitness function is arbitrarily set to $10^{20}$. 

Second when the CFL condition is violated, the  simulation is numerically unstable, and generates absurd values. In this case, the simulation is stopped, and the fitness function is arbitrarily set to a value larger than $10^6$. Note that a better solution would be to detect such violation before running the simulation, and to penalize the fitness by some amount that would be proportional to the actual violation. But it is numerically intractable to predict in advance if the CFL is going to be violated (see
\eq{eq:cfl}), and the numerical absurd values returned in case of numerical instability are not clearly correlated with the amount of violation either.\\

\paragraph*{\bf Initialization} 
The initial mean $\xmean{0}$ for CMA-ES is uniformly drawn in
$[-1,1]^{n}$, i.e., the parameters
$\KK_a$, $\NN^*_b$ and $\EE_c$ are uniformly drawn in the ranges given by the
expert. The initial step-size $\sigma_0$ is set to $0.3$.   Besides we
reject individuals of the population sampled outside the initial ranges. Unfeasible individuals are also rejected at initialization: at least one individual should be feasible to avoid random behavior of the
algorithm. In both cases, rejection is done by resampling until a ``good'' individual is got or a maximal number of sampling individuals is reached. Initial numbers of offspring $\lambda$ and parents $\mu$ are set to the default
values ($\lambda = \lfloor 4+3\log({\textrm{n}}) \rfloor$ and $\mu = \lfloor \lambda/2 \rfloor$).\\

\paragraph*{\bf Restarting and stopping criteria}
The algorithm stops if it reaches $5$ restarts, or a given fitness value
(typically a value between $10^{-9}$ and $10^{-15}$ for artificial problems, and adjusted for real data). 
Restart criteria (see Section \ref{CMA}) are RestartTolFun with
TolFun$=10^{-12} \times \sigma^{(0)}$, RestartTolX with
TolX$=10^{-12} \times \sigma^{(0)}$, RestartOnNoEffectAxis and
RestartCondCov with a limit upper bound of
$10^{14}$ for the condition number. The offspring size $\lambda$ is doubled after each restart and  $\mu$ is set equal to $\lfloor \lambda/2 \rfloor$.


\section{Results}
\label{sec:results}


\subsection{Validation using artificial data}
A first series of validation runs was carried out using simulated chromatograms. Each identification uses one or many experimental chromatograms. Because the same discretization is used for both the identification and the generation of the "experimental" data, one solution is known (the same isotherm that was used to generate the data), and the best possible fitness is thus zero.

\begin{figure}[Single component mixture]
 \centerline{\includegraphics[width=1.9in,angle=270]{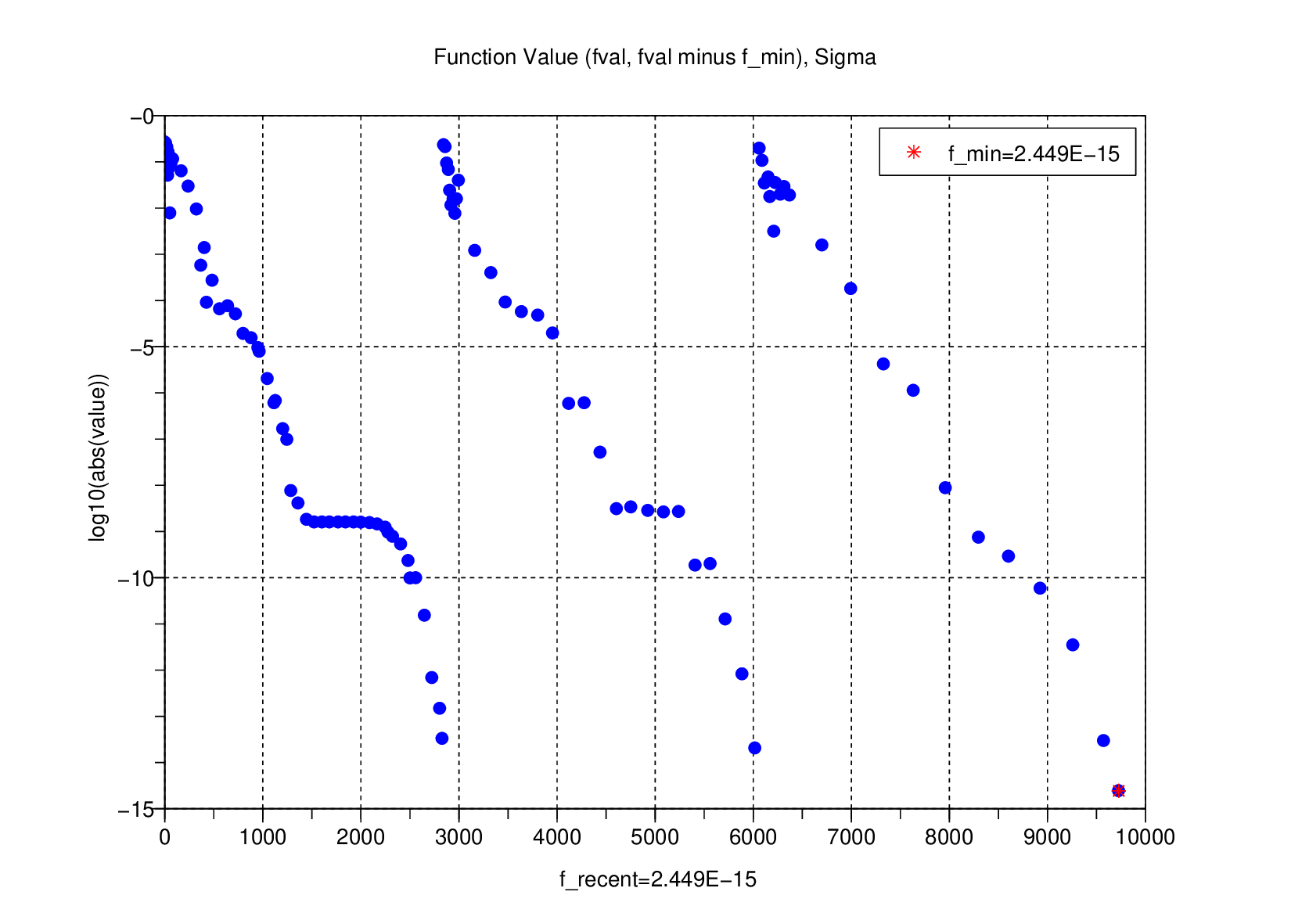}}
 \caption{Single component mixture, $1000$ time steps. Simulate a Lattice ($5$ parameters) and identify a
Lattice of degree $4$ ($5$ parameters): Best
fitness versus number of evaluations. The first run gave  a satisfactory solution but two restarts have been performed
 to reach a fitness value ($2.4 \ 10^{-15}$) lower than  $10^{-14}$.}
 \label{fig_sima}
 \end{figure}

\begin{figure}[Binary component mixture]
 \centerline{\includegraphics[width=1.9in,angle=270]{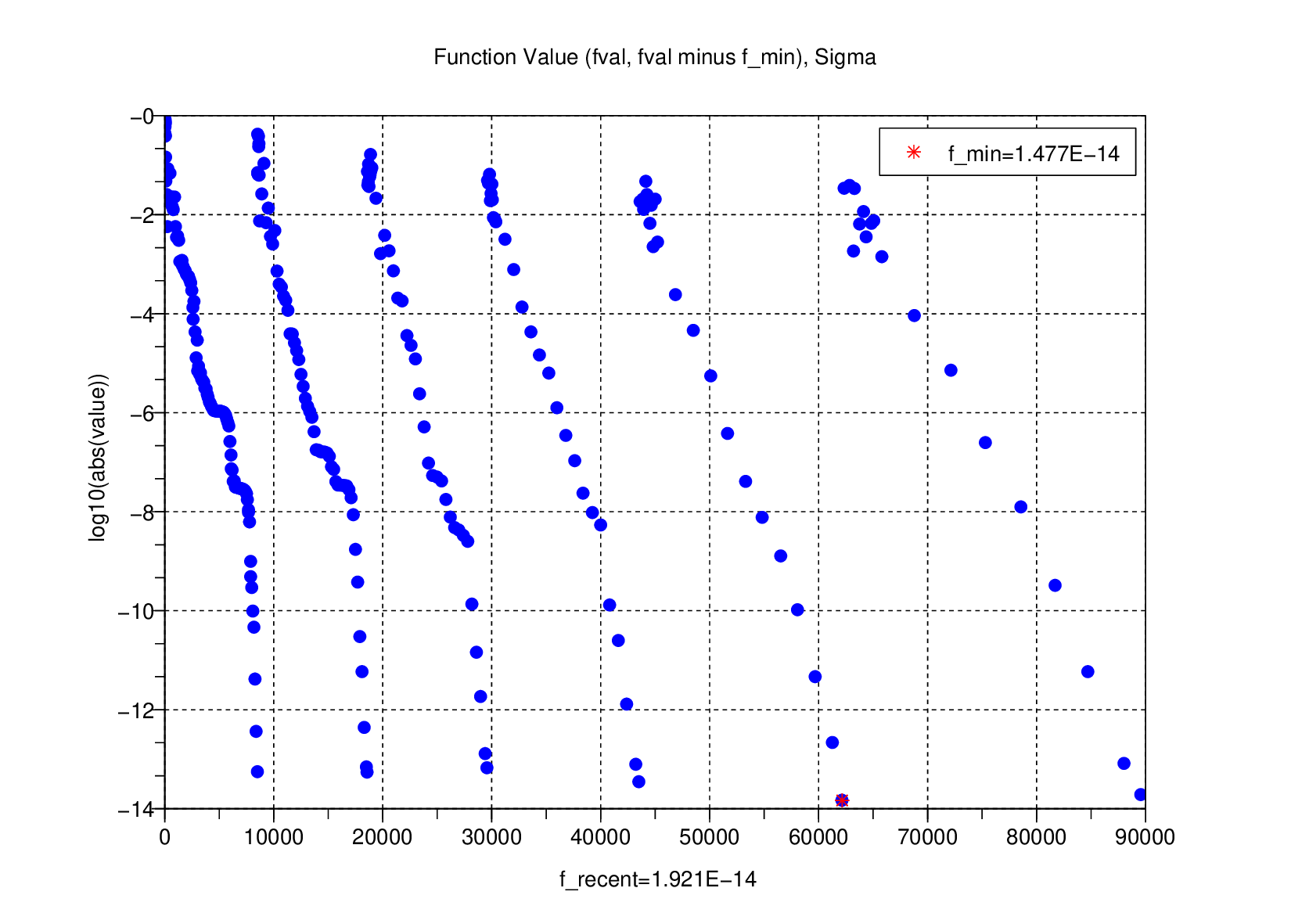}}
 \caption{Binary component mixture, $500$ time steps . Simulate a Langmuir ($3$ parameters) and identify a
Lattice of degree $3$ ($10$ parameters): Best
fitness versus number of evaluations.  The first run gave  a satisfactory solution
but the maximal number (here five) of  restarts have been performed attempting to reach a fitness value of 
$10^{-14}$, the best fitness value  ($1.4 \  10^{-14}$) 
was reached in the fourth restart. }
 \label{fig_simb}
 \end{figure}

\begin{figure}[Ternary component mixture]
\centerline{\includegraphics[width=1.9in,angle=270]{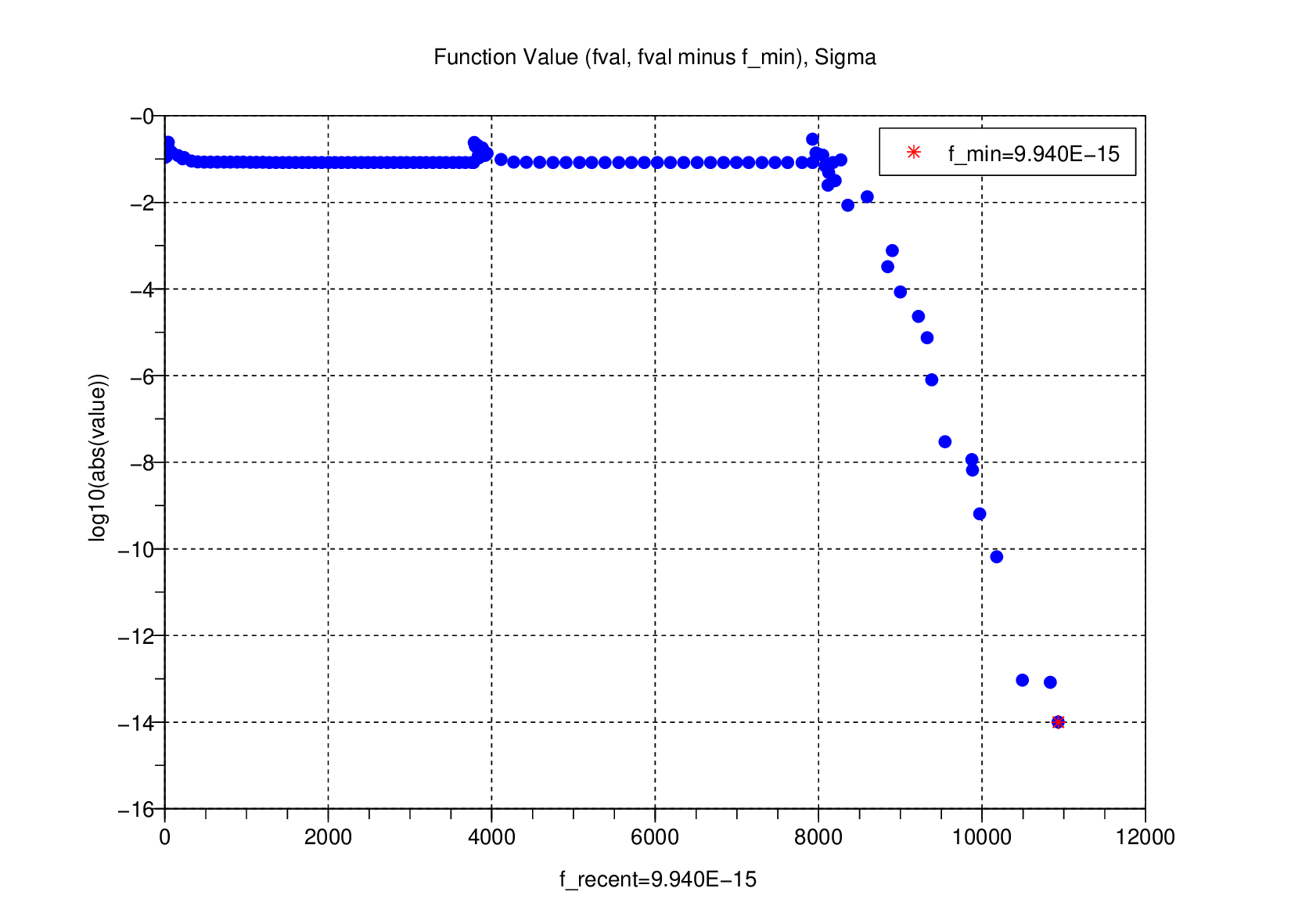}}
\caption{Ternary component mixture, $2000$ time steps. Simulate a Langmuir ($4$ parameters) and identify a
Langmuir ($4$ parameters): Best fitness versus
number of evaluations. Two restarts were necessary: Before the
second restart, CMA-ES is stuck in some local optima (fitness of order of $10^{-1}$),
in the second restart, the algorithm reaches a fitness value of $9.9 \ 10^{-15}$.}
\label{fig_simc}
\end{figure}

Several tests were run using different models for the isotherm,
different number of components, and different numbers of time
steps. In all cases, CMA-ES identified the correct parameters, {\it
i.e.} the fitness function reaches values very close to zero. In most cases, CMA-ES
did not need any restart to reach a precision of ($10^{-14}$), though this was necessary in a few cases. This happened when the whole population
remained unfeasible during several generations, or when the algorithm was stuck in
a local optimum. Figures~\ref{fig_sima},~\ref{fig_simb},~\ref{fig_simc} show typical evolutions during one run of the best fitness value with respect to the number of evaluations, for problems involving respectively 1, 2 or 3 components. Figure~\ref{fig_simc} is a case where restarting allowed the algorithm 
to escape a local optimum.

Specific tests were then run in order to study the influence of the expert guesses about both the ranges of the variables and the starting point of the algorithm possibly given by the chemical engineers: In CMA-ES, in a generation $g$,  offspring are drawn from a Gaussian distribution centered on the mean $\xmean{g}$. An expert guess for a good solution can hence be input as the mean of the first distribution $\xmean{0}$ that will be used to generate the offspring of the first generation. The results are presented in Table~\ref{table_comparALL}. First 3 lines give the probabilities that a given run converges (i.e., reaches a fitness value of $10^{-12}$), computed on 120 runs, and depending on the number of restarts (this probability of course increases with the number of restarts). The last line is the ratio between the average number of evaluations that were needed before convergence (averaged over the runs that did converge), and the probability of convergence: this ratio measures the performance of the different experimental settings, as discussed in details in \cite{Auger:2005b}.

The results displayed in Table~\ref{table_comparALL} clearly demonstrate that a good guess of the range of the variables is the most prominent factor of success: even without any hint about the starting point, all runs did reach the required precision without any restart. However, when no indication about the range is available, a good initial guess significantly improves the results, without reaching the perfect quality brought by tight bounds on the ranges: scaling is more important than rejecting unfeasible individuals at the beginning.\\

\begin{table}
\begin{center}
\renewcommand{\arraystretch}{1.3}
\caption{On the usefulness of Expert Knowledge: target values for Langmuir isotherm are here $(\KK_1,\NN^*) = (0.0388, 107)$. Expert range is $[0.01,0.05] \times [50,150]$, wide range is $[0.001,1] \times [50,150]$. The expert guess for the starting point is a better  initial mean (according to fitness value) than random. 
The first 3 lines give the  probabilities (computed over 120 runs) to reach a $10^{-12}$ fitness value within the given number of restarts. The last line is the ratio of the number of evaluations needed for convergence (averaged over the runs that did converge) by the probability of convergence after two restarts (line 3).}
\label{table_comparALL}
\begin{tabular}{c|c|c|c}
 Range &  Expert range  & Wide range & Wide range \\
 Starting point &  No guess  & No guess & Expert guess\\
\hline
no restart & $1$ & $0.84$ & $0.95$\\
$1$ restart & $1$  & $0.92$ & $0.97$\\
$2$ restarts & $1$ &  $0.95$ & $0.97$\\ 
\hline
\hline
Perf. & $601$ & $1015$ & $905$\\
\end{tabular}
\end{center}
\end{table}

\paragraph*{\bf Computational cost}
The duration of an evaluation depends on the discretization  of the numerical scheme (number of space- and time-steps), and on the number $n$  of unknown parameters to identify. Several runs were precisely timed to assess the dependency of the computational cost on both factors. The simple Langmuir isotherm was used to both generate the data and identify the isotherm. Only computational costs of single evaluations are reported, as the number of evaluations per identification heavily depends on many parameters, including the possible expert guesses, and in any case is a random variable of unknown distribution. All runs in this paper were performed on a 1.8GHz Pentium computer running with a recent Linux system.

For one component ($m=1$, $n=2$), and $100$, $500$ and $1000$ time steps, the  averages of the durations of a single evaluation are respectively $0.0097$, $0.22$, and $0.9$ seconds, fitting the theoretical quadratic increase with the number of time steps (though 3 sample points are too few to demonstrate anything!). This also holds for the number of space steps as the number of space steps is proportional  to the number of time steps due to the CFL condition. For an identification with a 1-component Langmuir isotherm, the total cost of the identification is on average $540$ seconds for a $1000$ time steps discretization.

When looking at the dependency of the computational cost on the number of unknown parameters, things are not that clear from a theoretical point of view, because the cost of each computation of the isotherm function also depends on the number of components and  on the number of experimental chromatograms to compare with. Experimentally, for, $2$, $3$ and $4$ variables, the costs of a single evaluation are respectively $0.9$, $1.04$, and $2.2$ seconds (for a $1000$ time steps discretization). For an identification, the total time is roughly  15 to 25 minutes for 2 variables, 40 to 60 minutes for 3 variables, and 1 to 2 hours for 4 variables.

\subsection{Experiments on real data}
The CMA-ES based approach has also been tested on a set of data taken from
\cite{cit:igor}. The mixture was composed of 3
chemical species: the benzylalcohol (BA), the 2-phenylethanol (PE) and
the 2-methylbenzylalcohol (MBA). Two real experiments have been performed with different proportions of injected mixtures, with respective proportions (1,3,1) and (3,1,0). Consequently, two real chromatograms have been provided.
For this identification, Qui\~{n}ones {\it et a.l.} \cite{cit:igor} have
used a {\it modified Langmuir} isotherm model in which each species
has a different saturation coefficient $\NN^*_i$:
\begin{equation}
\label{eq:lang-Igo}
\HH_i(c)=\frac{\NN_{i}^*}{1+\sum_{l=1}^3  \KK_l \  c_l }   \KK_i \  \cc_i, \ i=1, \ldots, 3.
\end{equation}
Six parameters are to be identified: $\NN_{i}^*$ and $\KK_i$,
for $i=1, \ldots, 3$. A change of variable has been made for those
tests so that the unknown parameters are in fact  $\NN_{i}^*$ and
$\KK_i^{'}$, where  $ \KK_i^{'}= \KK_i \, \NN_i$: those are the values that chemical engineers are able to experimentally measure.

\begin{figure}
\centerline{\subfigure[1:3:1, BA, PE]{\includegraphics[width=1.6in]{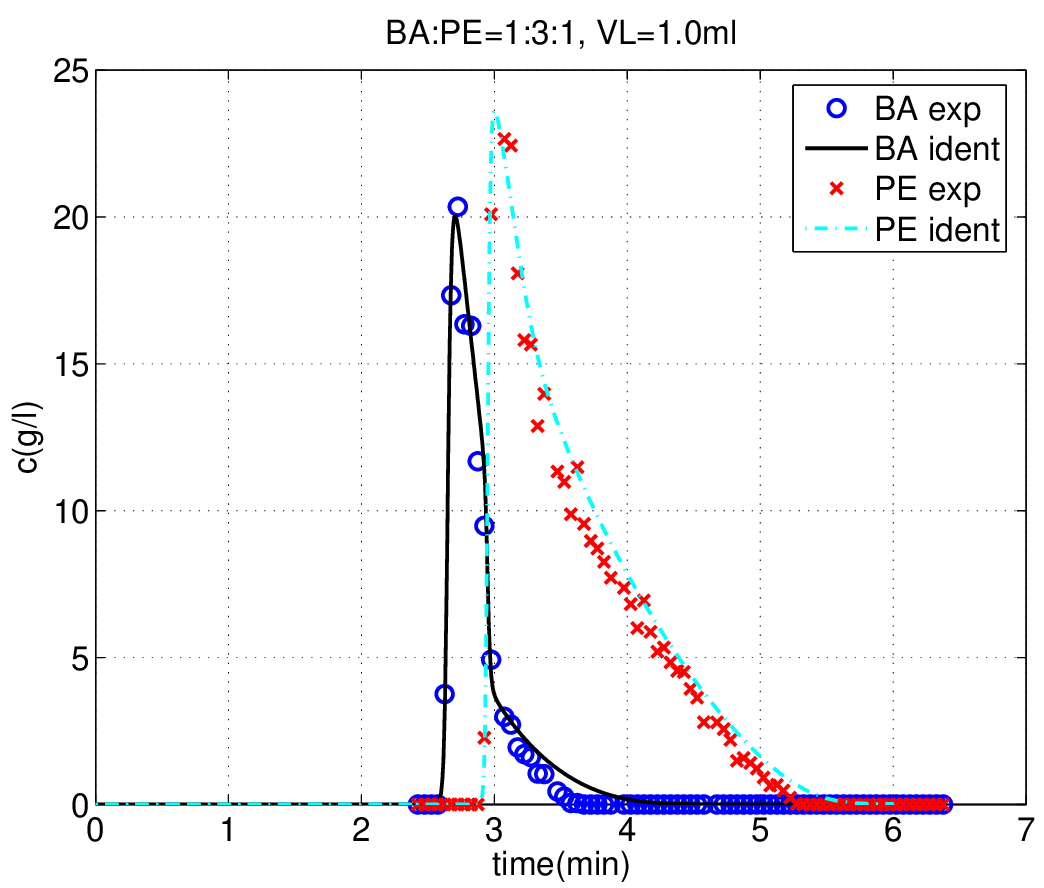}
\label{fig_sim61}}
\hfil
\subfigure[3:1, BA, PE]{\includegraphics[width=1.6in]{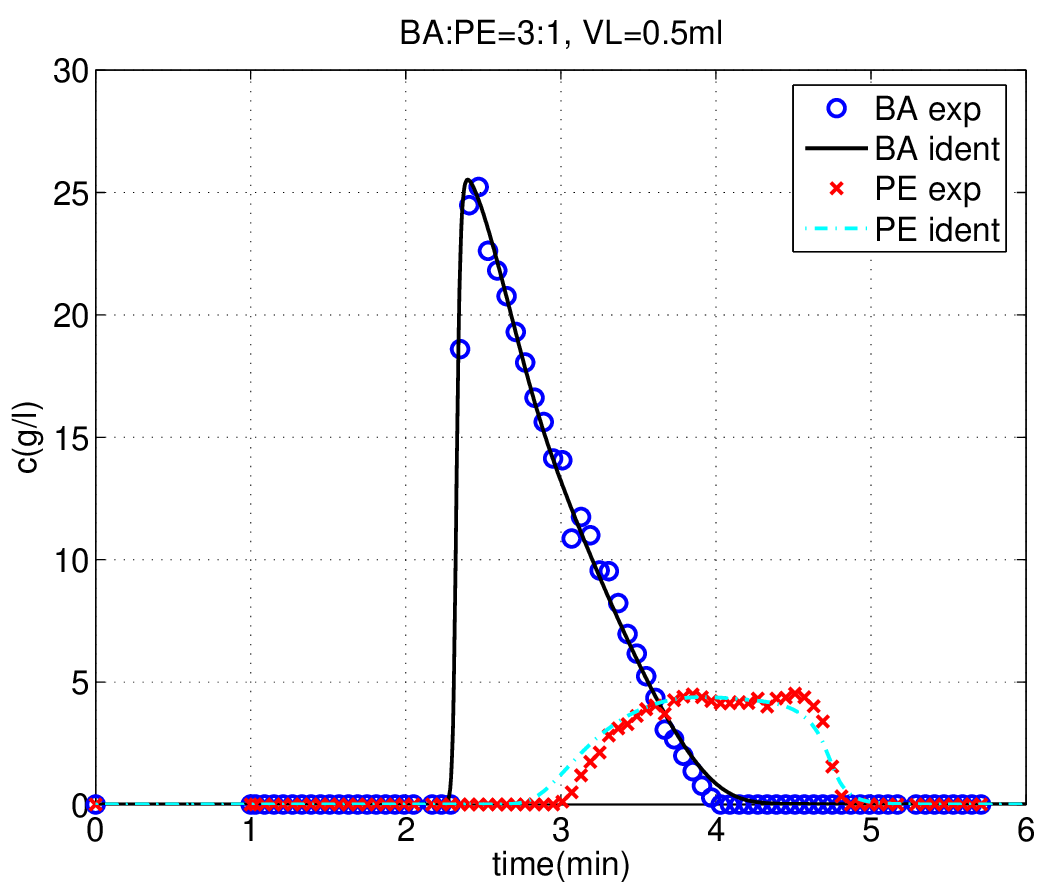}
\label{fig_sim62}}}
\centerline{\subfigure[1:3:1, MBA]{\includegraphics[width=1.6in]{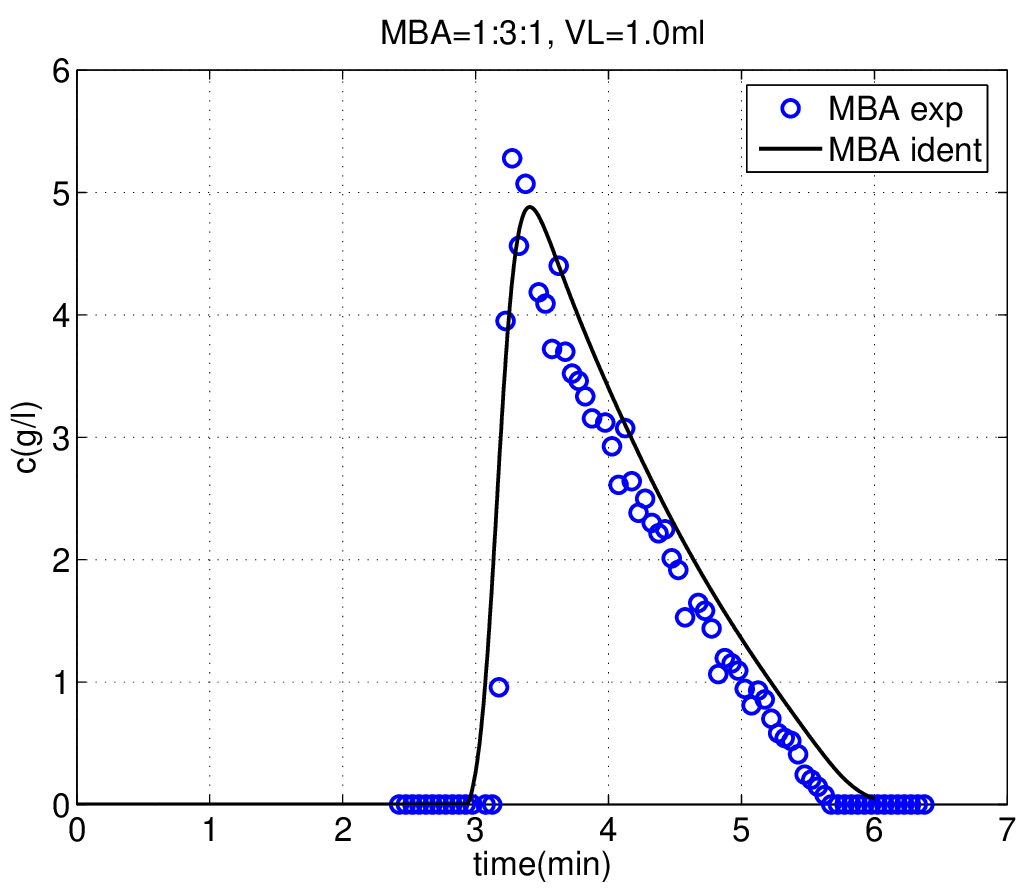}
\label{fig_sim63}}
\hfil
\subfigure[3:1, MBA]{\includegraphics[width=1.6in]{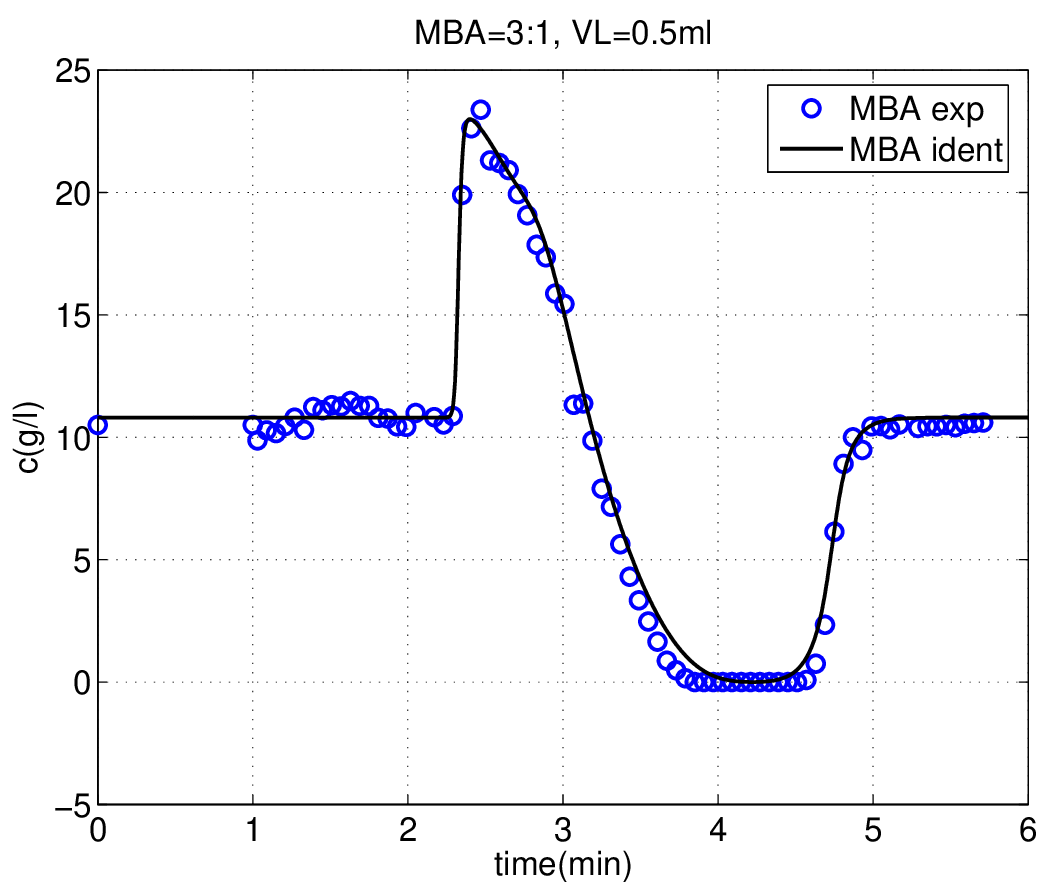}
\label{fig_sim64}}}
\caption{Experimental chromatograms (markers) and identified
chromatograms (continuous line) for the BA, BE and MBA species. Plots
on the left/right correspond to an injection with proportions
(1,3,1)/(3,1,0).}
\label{fig_sim6}
\end{figure}

Two series of numerical tests have been performed using a gradient-based method \cite{James:2007}: identification of the whole set of 6 parameters, and identification of the 3 saturation coefficients $\NN_{i}^*$ only, after setting the Langmuir coefficients  to the experimentally  measured values 
$(\KK_1^{'},\KK_2^{'},\KK_3^{'}) = (1.833, 3.108, 3.511)$.
The initial ranges
used for CMA-ES are $[60,250]\times[60,250] \times [60,250]$
(resp. $[1.5,2.5]\times[2.7,3.7] \times [3,4] \times [90,200] \times
[100,200] \times [100,210]$) when optimizing 3 parameters (resp. 6 parameters). Comparisons
between the two experimental chromatograms and those resulting from CMA-ES
identification for the two experiments
are shown in Figure~\ref{fig_sim6}, for the $6$-parameters case. The corresponding plots in the $3$-parameters case are visually identical though the fitness value is slightly lower than in the  $6$-parameters case (see Tables \ref{table_compar_isoth} and \ref{table_compar_isoth2}).
But another point of view on the
results is given by the comparison between the identified isotherms
and the (few) experimental values gathered by the chemical engineers. The
usual way to present those isotherms in chemical publications is that
of Figure~\ref {fig_isoth6}: the absorbed quantity $\HH(\cc)_i$ of
each component $i=1,2,3$ is displayed as a function of the total
amount of mixture $(\cc_1+\cc_2+\cc_3)$, for five different
compositions of the mixture \cite{James:2007}.  Identified
(resp. experimental) isotherms are plotted in Figure~\ref {fig_isoth6}
using continuous lines (resp. discrete markers), for the $6$-parameters case. Here again the corresponding plots for the $3$-parameters case are visually identical.


%




\begin{table}
\begin{center}
\renewcommand{\arraystretch}{1.3}
\caption{Comparing CMA-ES and gradient: the $3$-parameters case.
Solution ( line 1) and
associated fitness values ( line 2) for the modified Langmuir model
(\eq{eq:lang-Igo}).  Line 3: For CMA-ES, "median (minimal)" number of fitness evaluations (out of 12 runs) needed to reach the corresonding fitness value on  line 2. For gradient, "number of fitness evaluations -- number of gradient evaluations" for the best
of the 10 runs described in \cite{James:2007}.
}
\label{table_compar_isoth}
\begin{tabular}{@{$\,$}c@{$\,$}|cc@{$\,$}c@{$\,$}}
 &  \multicolumn{2}{c}{CMA-ES}   & Gradient  \\
\hline
$\NN_i^{*}$ & \multicolumn{2}{c}{(120.951,135.319,165.593)} & (123.373,135.704,159.637) \\
Fitness $\times$ $10^{3}$ & 8.96 & {\bf 8.78}  & 8.96  \\
\! \# \! Fit evals. \! & 175 (70) & 280 (203)   &  140 -- 21 \\
\end{tabular}
\end{center}
\end{table}

\subsection{Comparison with a Gradient Method}
CMA-ES results have then been compared with those of the gradient method from
\cite{James:2007}, using the same data case of ternary mixture taken
from \cite{cit:igor} and described in previous Section. Chromatograms
found by CMA-ES are, according to the fitness (see
Tables~\ref{table_compar_isoth} and~\ref{table_compar_isoth2}), closer
to the experimental ones than those obtained with the gradient method.
Moreover, contrary to the gradient algorithm, all 12 independent runs of CMA-ES converged to the same point. Thus, no variance is to be reported on Tables~\ref{table_compar_isoth} and~\ref{table_compar_isoth2}. Furthermore, there
seems to be no need, when using CMA-ES, to fix the $3$ Langmuir
coefficients in order to find good results: when optimizing all $6$
parameters, the gradient approach could not reach a value smaller than
$0.01$, whereas the best fitness found by CMA-ES in the same context
is $8.32 \, 10^{-3}$ (Table \ref{table_compar_isoth2}).

\begin{figure}
\centerline{\subfigure[BA]{\includegraphics[width=1.63in]{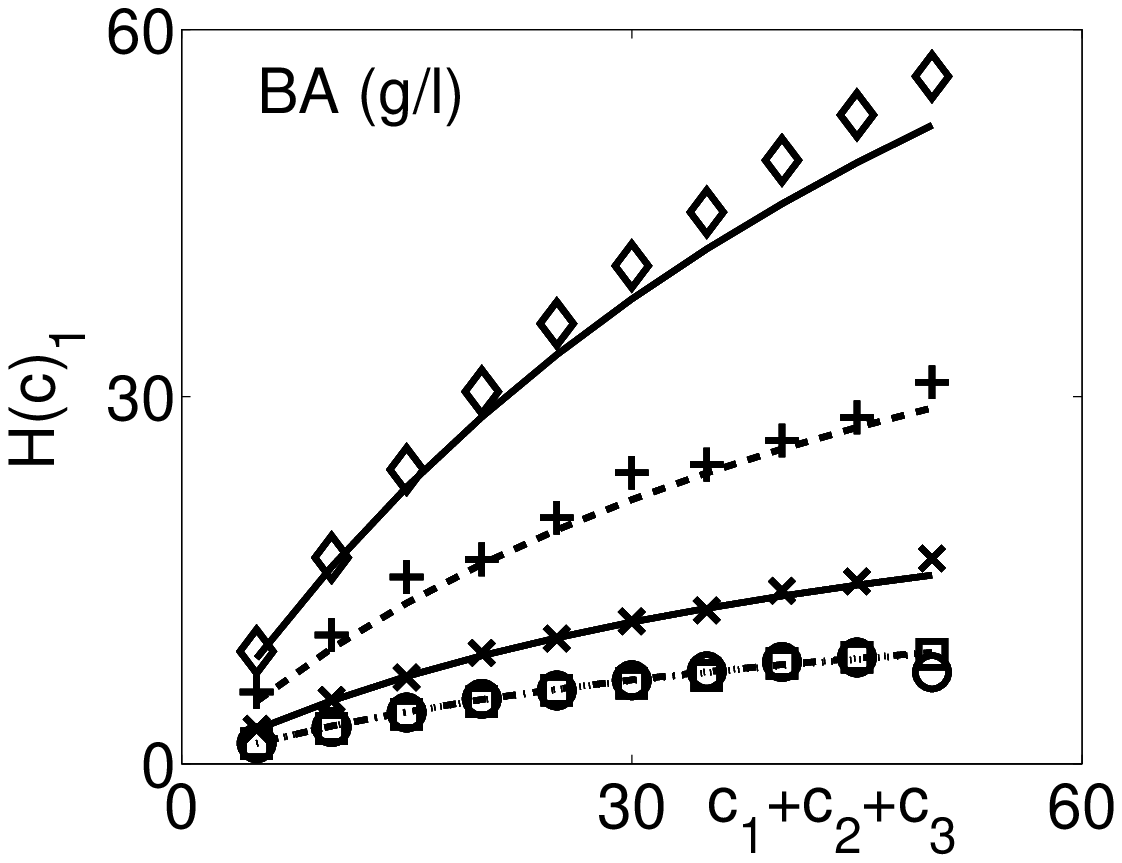}
\label{fig_isoth61}}
\hfil
\subfigure[MBA]{\includegraphics[width=1.63in]{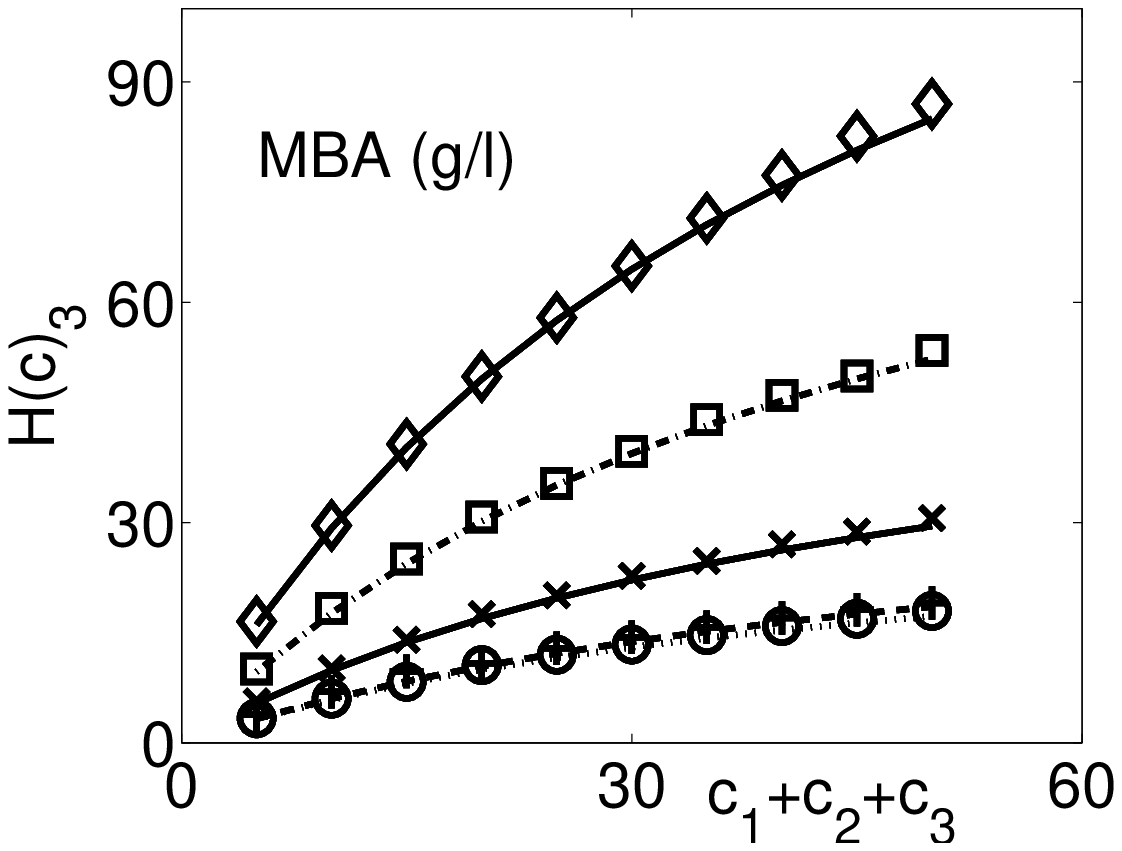}
\label{fig_isoth62}}}
\centerline{\subfigure[PE]{\includegraphics[width=1.63in]{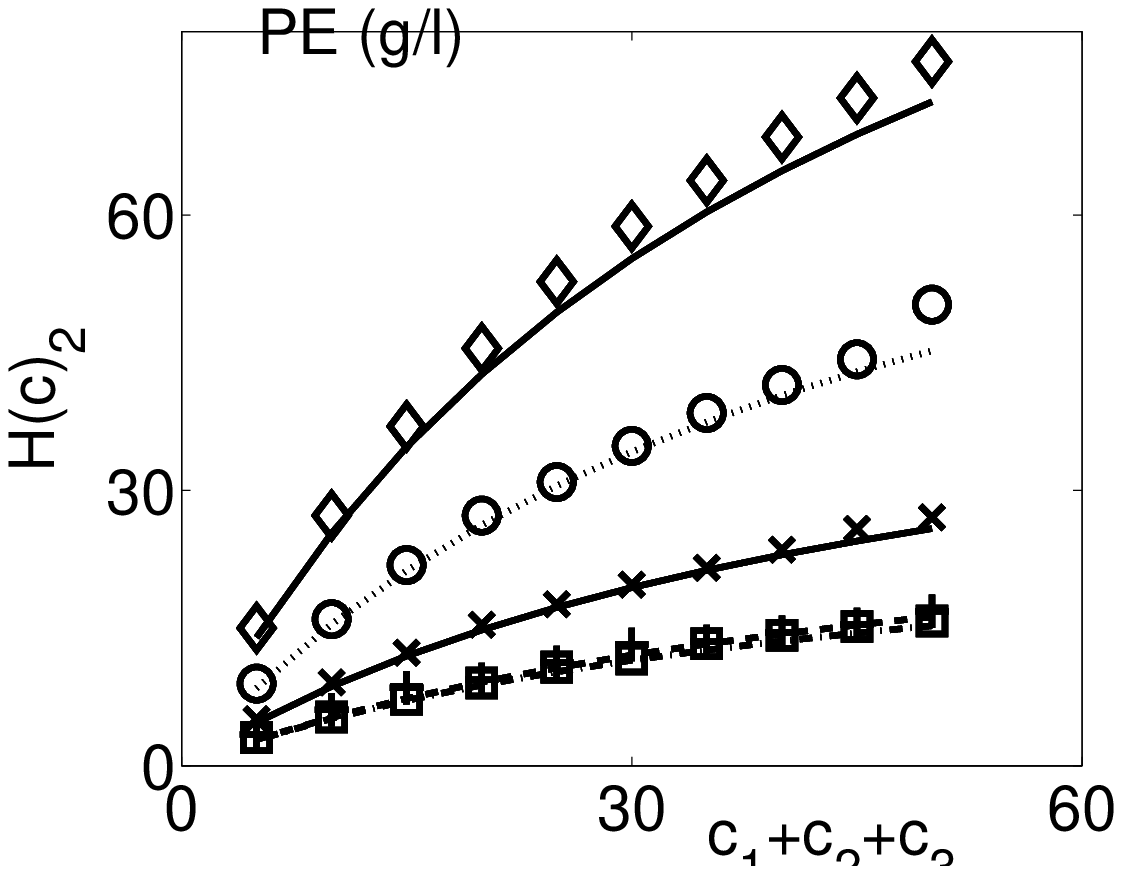}
\label{fig_isoth63}}}
\caption{Isotherms associated to parameters values of Table~\ref{table_compar_isoth2} (continuous line) and experimental ones (markers) versus total amount of the mixture for different proportions of the component in the injected concentration \cite{James:2007}.}
\label{fig_isoth6}
\end{figure}

\begin{table}
\begin{center}
\renewcommand{\arraystretch}{1.3}
\caption{Comparing CMA-ES and gradient: the $6$-parameters case. 
Solutions ( lines 1 and 2) and
associated fitness values ( line 3) for the modified Langmuir model
(\eq{eq:lang-Igo}).
}
\label{table_compar_isoth2}
\begin{tabular}{c|cc}
 &  CMA-ES & Gradient \\
\hline
  $K_i^{'}$ & (1.861,3.120,3.563) & (1.780,3.009,3.470)  \\ 
  $N_i^{*}$ & (118.732,134.860,162.498) & (129.986,141.07,168.495) \\
  Fitness  $\times$ $10^{3}$ & {\bf 8.32}  & 10.7  \\
\end{tabular}
\end{center}
\end{table}

%

Finally, when comparing the identified isotherms to the experimental
ones (figure \ref{fig_isoth6}), the fit is clearly not very satisfying
(similar deceptive results were obtained with the gradient method in
\cite{James:2007}): Fitting both the isotherms and the chromatograms
seem to be contradictory objectives. Two directions can lead to some
improvements in this respect: modify the cost function $\cost$ in
order to take into account some least-square error on the isotherm as
well as on the chromatograms; or use a multi-objective approach. Both
modifications are easy to implement using Evolutionary Algorithms (a
multi-objective version of CMA-ES was recently proposed
\cite{CMAES-Multi}), while there are beyond what gradient-based
methods can tackle. However, it might also be a  sign that the
modified Langmuir model that has been suggested for the isotherm
function is not the correct one. \\

\paragraph*{\bf Comparison of convergence speeds}

Tables~\ref{table_compar_isoth} and~\ref{table_compar_isoth2} also give an idea of the respective computational costs of both methods on the same real data.
For the best run out of 10, the
gradient algorithm reached its best fitness value after $21$
iterations, requiring on average $7$ evaluations per iteration 
for the embedded line search. Moreover, the computation of the gradient itself is costly -- roughly estimated to 4 times that of the fitness function.
Hence, the total cost of the gradient algorithm can be considered to be larger than $220$ fitness evaluations. To reach the same fitness value ($8.96 \, 10^{-3}$), CMA-ES only needed 175 fitness evaluations (median value out of 12 runs). To converge to its best value ($8.78 \, 10^{-3}$, best run out of 12) CMA-ES needed 280 fitness evaluations. Those results show that the best run of  the gradient algorithms needs roughly the same amount of functions evaluations than CMA-ES to converge. Regarding the robustness issue, note that CMA-ES always reached the same fitness value, while the 10 different runs of the gradient algorithm from 10 different starting points gave 10 different solutions: in order to assess the quality of the solution, more runs are needed for the gradient method than for CMA-ES!

\section{Conclusions}
This paper has introduced the use of CMA-ES for the parametric identification of isotherm functions in chromatography. Validation tests on simulated data were useful to adjust the (few) CMA-ES parameters, but also demonstrated the
importance of expert knowledge: choice of the type of isotherm, ranges for the different parameters, and possibly some initial guess of a not-so-bad solution.

The proposed approach was also applied on real data and compared to previous work
using gradient methods. On this data set, the best fitness found by
CMA-ES is better than that found by the gradient approach. Moreover, the results obtained with CMA-ES are far more
robust: (1) CMA-ES always converges to the same values of the isotherm parameters, independently of its starting point; (2) CMA-ES can handle the full problem that the gradient method failed to efficiently  solve: there is no need when using CMA-ES to use experimental values of the Langmuir parameters in order to obtain a satisfactory fitness value.
Note that the fitness function only takes into account  the fit of the chromatograms, resulting in a poor fit on the isotherms. The results confirm the
ones obtained with a gradient approach, and suggest to either incorporate some measure of isotherm fit in the fitness, or to try some multi-objective method -- probably the best way to go, as both objectives (chromatogram and isotherm fits) seem somehow contradictory.




\section*{Acknowledgments}

This work was supported in part by MESR-CNRS ACI NIM Chromalgema. The
authors would like to thank Nikolaus Hansen for the
Scilab version of CMA-ES, and for his numerous useful comments.

\def\V{\rm vol.~}
\def\N{no.~}
\def\pp{pp.~}
\def\Pot{\it Proc. }
\def\IJCNN{\it International Joint Conference on Neural Networks\rm }
\def\ACC{\it American Control Conference\rm }
\def\SMC{\it IEEE Trans. Systems\rm, \it Man\rm, and \it Cybernetics\rm }

\def\handb{ \it Handbook of Intelligent Control: Neural\rm, \it
    Fuzzy\rm, \it and Adaptive Approaches \rm }



\bibliographystyle{plain}
\bibliography{cec07}

\end{document}